\documentclass[12pt]{article}
\usepackage{amsmath,amssymb,amsthm, tikz}
\usepackage{amsfonts, amsmath}
\usepackage[english]{babel}
\newtheorem{theorem}{Theorem}[section]

\theoremstyle{remark}

\numberwithin{equation}{section}

\newcommand{\restr}{\lfloor}
\newcommand\eps\varepsilon

\newcommand{\R}{\mathbb R}

\newcommand{\beq}{\begin{equation}}
\newcommand{\eeq}{\end{equation}}

\begin{document}
\title{\textbf{Some recent results on the Dirichlet problem for $(p,q)$-Laplace equations}}
\author{
\bf Salvatore A. Marano\thanks{Corresponding author}, Sunra
J.N. Mosconi\\
\small{Dipartimento di Matematica e Informatica,
Universit\`a degli Studi di Catania,}\\
\small{Viale A. Doria 6, 95125 Catania, Italy}\\
\small{\it E-mail: marano@dmi.unict.it, mosconi@dmi.unict.it}\\
}
\date{}
\maketitle
\begin{abstract}
A short account of recent existence and multiplicity theorems on the Dirichlet problem for an elliptic equation with $(p,q)$-Laplacian in a bounded domain is performed. Both eigenvalue problems and different types of perturbation terms are briefly discussed. Special attention is paid to possibly coercive, resonant, subcritical, critical, or asymmetric right-hand sides.
\end{abstract}
\vspace{2ex}
\vspace{2ex}

\section{Introduction}
Let $\Omega$ be a bounded domain in $\R^N$ with a $C^2$-boundary $\partial\Omega$, let $1<q\leq p<N$, and let $\mu\in\R^+_0$. Consider the Dirichlet problem
\begin{equation}\label{problema}
\left\{
\begin{array}{ll}
-\Delta_p u-\mu \Delta_q u=g(x, u)\quad & \mbox{in}\quad\Omega,\\
u=0\quad & \mbox{on}\quad\partial\Omega,\\
\end{array}
\right.
\end{equation}
where $\Delta_r$, $r>1$, denotes the $r$-Laplacian, namely
$$\Delta_r u:={\rm div}(|\nabla u|^{r-2}\nabla u)\quad\forall\, u\in W^{1,r}_0(\Omega),$$
$p=q$ iff $\mu=0$, while $g:\Omega\times\R\to\R$ satisfies Carath\'{e}odory's conditions. The non-homogeneous differential operator $Au:=\Delta_p u+\mu\Delta _q u$ that appears in \eqref{problema} is usually called $(p,q)$-Laplacian. It stems from a wide range of important applications, including biophysics \cite{Fi}, plasma physics \cite{Wil}, reaction-diffusion equations \cite{Ar,CI}, as well as models of elementary particles \cite{BDFP,BFP,BMV, De}. That's why the relevant literature looks daily increasing and numerous meaningful papers on this subject are by now available.

This survey provides a short account of some recent existence and multiplicity results involving \eqref{problema}. To increase readability, we chose to report special but significant cases of more general theorems, and our statements are often straightforward extensions of known results. We refer to the original papers for the most up-to-date research on this topic.

Section 2 contains  basic properties of the operator $u\mapsto -\Delta_p u-\mu\Delta_q u$, $u\in W^{1,p}_0(\Omega)$. In particular, a scaling argument shows that if $\mu>0$ then there is no loss of generality in assuming $\mu:=1$, which we will make throughout the work. Eigenvalue problems are treated in Section 3, where, accordingly,
$$g(x,t):=\alpha|t|^{p-2}t+\beta|t|^{q-2}t\quad\forall\, (x,t)\in\Omega\times\R,$$
with $(\alpha,\beta)\in\R^2$. The results we present are chiefly taken from \cite{BT,CDG,MaPa,T0,T} and concern bounded domains; for the whole space see, e.g., \cite{BZ0,BZ,HeLi}. Section 4 collects existence and multiplicity theorems, sometimes with a precise sign information. We shall suppose $g$ purely autonomous as well as of the type
$$g(x,t):=\alpha|t|^{p-2}t+\beta|t|^{q-2}t+f(t),\quad (x,t)\in\Omega\times\R,$$
where the perturbation $f$ lies in $C^1(\R)$ and exhibits a suitable growth rate at $\pm\infty$ and/or at zero. Results are diversified according to the asymptotic behavior at infinity, frequently under the simplified condition
$$-\infty<\lim_{|t|\to +\infty}\frac{f(t)}{|t|^{r-2}t}<+\infty$$
for some $r>1$. In particular, Section 4.1 concerns the $(p-1)$-sublinear case $r\le p$, possibly with coercivity and/or resonance additional assumptions. Section 4.2 treats the $(p-1)$-superlinear case $r>p$ and both subcritical, i.e.,  $r<p^*:=\frac{Np}{N-p}$, or critical, namely $r=p^*$, situations are examined. Finally, Section 4.3 deals with asymmetric nonlinearities, meaning that the asymptotic behavior at $-\infty$ and $+\infty$ is different. Due to limited space, we relied only on \cite{CMP,CDG,FQ}, \cite{MMP}--\cite{MePe}, \cite{MuPa}--\cite{PaRa3}, \cite{PRR,PaWi2}, \cite{YangBai}--\cite{YY}, where bounded domains are considered, and refer to \cite{BZ,CEM,F,HeLi} for the case $\Omega:=\R^N$. 

We surely forgot to mention here significant works, something of which we apologize in advance.
\section{Basic properties of the $(p, q)$-Laplacian}
Let $1<q\le p< N$ and let $\mu\in\R^+_0$. We denote by $p'$ the dual exponent of $p$, i.e., $p':=p/(p-1)$, while $p^*$ is the Sobolev dual in dimension $N$, namely 
$$p^*:=\frac{Np}{N-p}.$$
If $r\in [1,+\infty)$ then
\[
\|u\|_r:=\left(\int_\Omega|\nabla u|^r\, dx\right)^{1/r}\quad\forall\, u\in W^{1,r}_0(\Omega),\quad  |u|_r:=\left(\int_\Omega|u|^r\, dx\right)^{1/p}\quad \forall\, u\in L^r(\Omega).
\]
Write $A_{p, q}^\mu$ for the differential of the $C^1$, strictly convex functional
\[
u\mapsto \frac{1}{p}\|u\|^p_p+\frac{\mu}{q}\|u\|_q^q, \quad u\in W^{1,p}_0(\Omega).
\]
Hence,
 \[
 \langle A_{p, q}^\mu(u), v\rangle=\int_\Omega\left(|\nabla u|^{p-2}\nabla u\cdot\nabla v+\mu|\nabla u|^{q-2}\nabla u\cdot\nabla v\right) dx\quad\forall\, u,v\in W^{1,p}_0(\Omega),
 \]
where, as usual, $\langle\cdot,\cdot\rangle$ indicates the duality coupling between $W^{1,p}_0(\Omega)$ and its topological dual $W^{-1,p'}(\Omega)$. Setting $\mu:=0$ we obtain the $p$-Laplacian $A_p$ while $\mu=1$ produces the so-called $(p,q)$-Laplacian. Let us now discuss some features of $A_{p,q}^\mu$.\\
 
\noindent {\em - Scaling}\\
The case $q=2$, $\mu\to 0^+$ naturally gives an elliptic regularization procedure for solutions $\bar u\in W^{1,p}_0(\Omega)$ of the equation $A_{p}(u)=g$, with $g\in W^{-1,p'}(\Omega)$, in the sense that solutions $u_\mu$ to $A_{p, 2}^\mu(u)=g$ usually posses better regularity properties for $\mu>0$ than for $\mu=0$ but, at the same time, $u_\mu\to\bar u$ strongly as $\mu\to 0^+$.\\
A scaling argument shows that the solutions of $A_{p, q}^\mu(u)=g$ actually solve an equation of the type $A_{p, q}^1(v)=h$, with explicit $v$ and $h$. So, the perturbation parameter $\mu$ can be avoided. To see this, for every $u\in W^{1,p}_0(\Omega)$ and $\lambda>0$ we set
$$(u)_\lambda(x):=u(\lambda x),\quad x\in\lambda^{-1}\Omega,$$
where $\lambda^{-1}\Omega:=\{\lambda^{-1}x:x\in \Omega\}$. Evidently, $(u)_\lambda\in W^{1,p}_0(\lambda^{-1}\Omega)$. By duality, given any $g\in W^{-1, p'}(\lambda^{-1}\Omega)$, one defines
 \[
 \langle (g)_\lambda, v\rangle:=\langle g, (v)_{\lambda^{-1}}\rangle\quad \forall\, v\in W^{1,p}_0(\lambda^{-1}\Omega),
 \]
whence $(g)_\lambda\in W^{-1, p'}(\lambda^{-1}\Omega)$. Now, if $A_{p, q}^\mu(u)=g$ for some $u\in W^{1,p}_0(\Omega)$ and $g\in W^{-1, p'}(\Omega)$ then, through a direct computation changing variables,
 \[
 A_{p, q}^1((u)_\lambda)=\lambda^{p}\big(A_{p}(u)\big)_{\lambda}+\lambda^q \big(A_q(u)\big)_{\lambda}=\lambda^p \big(A_{p, q}^{\lambda^{q-p}}(u)\big)_{\lambda},
 \]
so that, letting $\lambda:=\mu^{1/(q-p)}$, it holds
\[
A_{p, q}^\mu(u)=g\quad \Leftrightarrow\quad A_{p, q}^1((u)_{\mu^{1/(q-p)}})=\mu^{p/(q-p)}(g)_{\mu^{1/(q-p)}}.
\]
Therefore, henceforth,  \textit{we shall pick $\mu:=1$ and, to simplify notation, put $A_{p, q}:=A_{p, q}^1$.}\\

\noindent {\em - Homogeneity}\\
The operator $A_{p, q}^\mu$ is not homogeneous whenever $p\neq q$. This prevents to apply, in a simple way, available tools from critical point theory. However, a reasoning analogous to the one exploited before allows to take advantage from such a trouble. Indeed, if $u\in W^{1,p}(\Omega)$ is a solution of the quasilinear equation $A_{p, q}^\mu(u)=g(\cdot, u)$ then 
\[
A_{p, q}(ku)=k^{p-1}A_p(u)+k^{q-1}A_q(u)=k^{p-1}A_{p, q}^{k^{q-p}}(u).
\]
Setting $\mu:=k^{q-p}$, the function $v:=\mu^{1/(q-p)}u$ solves
\[
A_{p, q}(v)=\mu^{\frac{p-1}{q-p}}g(\cdot,u)=g_{\mu}(\cdot, v), \quad \text{with}\quad g_\mu(x, t):=\mu^{\frac{p-1}{q-p}}g(x, \mu^{\frac{1}{p-q}}t).
\]
We thus achieve an equation for $A_{p, q}$ where the reaction exhibits a different asymptotic behavior both at zero and at infinity.\\ 

\noindent {\em{- Eigenvalues of $A_p$}} \\
Eigenvalues of $A_p$ play a key r\^{o}le in solving quasilinear equations where $A_{p,q}$ appears. So, some basic properties will be recalled. Denote by $\sigma (A_p, W^{1,p}_0(\Omega))$ the spectrum of $A_p$ with (zero) Dirichlet boundary conditions. A whole sequence of variational eigenvalues
$$\sigma_p:=\{\lambda_k(p)\}\subseteq \sigma (A_p, W^{1,p}_0(\Omega))$$
can be constructed via the Ljusternik-Schnirelmann minimax scheme with $\mathbb{Z}_2$ co-homological index \cite[Section 4.2]{PAR}. If $p=2$ it reduces to the usual spectrum of $(-\Delta,H^1_0(\Omega))$, while for general $p>1$ one has
\[
0<\lambda_1(p)=\inf_{u\in W^{1,p}_0(\Omega)\setminus \{0\}}\frac{\|u\|_p^p}{|u|_p^p}<\lambda_2(p)=
\inf\{\lambda\in \sigma_p:\,\lambda\neq\lambda_1(p)\}.
\]
Moreover,  as long as $\Omega$ is connected, the eigenspace corresponding to $\lambda_1(p)$ is one-dimensional and  spanned by a positive function $u_{1,p}\in C^{1,\alpha}(\Omega)$, whereas eigenfunctions relative to  higher eigenvalues must be nodal.\\

\noindent {\em{- Regularity}}\\
The regularity theory for $A_{p, q}$ is well established after the works of Lieberman in the eighties (see, e.g., \cite{L}) and parallels that concerning $A_p$. In particular, weak solutions $u\in W^{1,p}_0(\Omega)$ to the equation
\begin{equation}\label{eq}
A_{p, q}(u)=g(\cdot, u),
\end{equation}
for some Carath\'{e}odory function $g:\Omega\times \R\to \R$ obeying 
\begin{equation}\label{subcritical}
|g(x, t)|\le C(1+|t|^{p^*-1}),
\end{equation}
are $C^{1,\alpha}$ up to the boundary. Under natural conditions on $g$, a strong maximum principle as well as the Hopf boundary point lemma hold true. Precisely, define
\[
C^1_0(\overline{\Omega}):=\{u\in C^1(\overline{\Omega}):\;u\restr_{\partial\Omega}=0\},
\]
Its positive cone
\[
C_+:=\{u\in C^1_0(\overline{\Omega}):\; u(x)\geq 0
\text{ in $\overline{\Omega}$}\}
\]
has a nonempty interior given by
\[
{\rm int } (C_+)=\left\{u\in C_+:\; u(x)>0\;\;\forall\, x\in\Omega,\;\;\frac{\partial u}{\partial n}(x)<0\;\;\forall\, x\in\partial\Omega\right\},
\]
where $n(x)$ denotes the outward unit normal to $\partial\Omega$ at $x$. Because of \cite[Theorems 5.5.1 and 5.3.1]{PS}, if there exists $\delta>0$, $C>0$ such that
\[
 tg(x, t)\ge - Ct^q\quad\mbox{in}\quad\Omega\times[-\delta,\delta]
\]
then, under \eqref{subcritical}, nonnegative nontrivial solutions to \eqref{eq} actually lie in ${\rm int } (C_+)$. A similar statement holds true for nonpositive solutions.
\section{Eigenvalue problems}\label{Sect3}

In analogy with the Fu\v cik spectrum, the eigenvalue problem for $A_{p, q}$ with homogeneous Dirichlet boundary conditions on a bounded domain $\Omega$ consists in finding all $(\alpha, \beta)\in \R^2$ such that the equation
\begin{equation}\label{eig}
A_{p, q}(u)=\alpha |u|^{p-2}u+\beta |u|^{q-2}u
\end{equation}
possesses a nontrivial weak solution $u\in  W^{1,p}_0(\Omega)$. The set of such $(\alpha, \beta)$ is called the $(p, q)$-spectrum of $A_{p, q}$ and denoted by $\sigma_{p, q}$. Additionally, $\sigma_{p, q}^+$ indicates the set of $(\alpha,\beta)\in\sigma_{p,q}$ for which there exists a {\em positive} solution to \eqref{eig}. This problem can evidently be recasted in the more general framework of weighted eigenvalues, namely to the equation
$$A_{p, q}(u)=\xi |u|^{p-2}u+\eta |u|^{q-2}u,$$
where $\xi$ and $\eta$ are bounded functions, and most of the results presented here have suitable weighted (sometimes even sign-changing) variants. For the sake of simplicity, we will focus on the constant coefficient case. Nevertheless, a full description of
$\sigma_{p, q}$ is out of reach, since it clearly presents additional difficulties with respect the comparatively simpler case of the spectrum of $A_p$, a well-known open problem until today. 

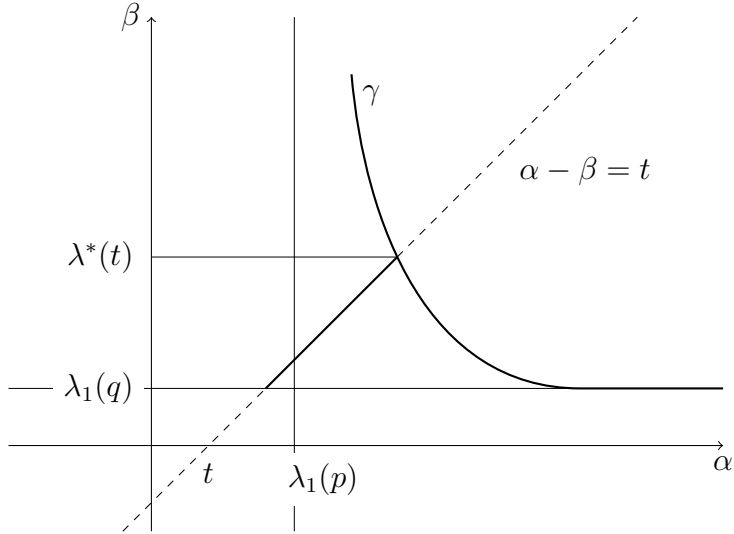
\begin{figure}
\centering
\begin{tikzpicture}[scale=1.9]
\draw[->] (-1, 0) -- (4, 0);
\draw (4, 0) node[below]{$\alpha$};
\draw[->] (0, -0.6) -- (0, 3);
\draw (0, 3) node[left]{$\beta$};
\draw[very thin] (-1, 0.4) -- (3, 0.4);
\draw[very thin] (1, -0.6) -- (1,3);
\draw (1.2,-0.05) node[below, fill=white]{$\lambda_1(p)$};
\draw (-0.05, 0.4) node[left, fill=white]{$\lambda_1(q)$};
\draw (0.4, -0.05) node[below]{$t$};
\draw[thin, dashed] (-0.2,-0.6) -- (0.8 ,0.4);
\draw[thick] (0.8, 0.4) -- (1.72, 1.32);
\draw[thin, dashed] (1.72, 1.32) -- (3.4, 3);
\draw (2.5, 2.1) node[below right]{$\alpha-\beta=t$};
\draw[thick] (1.4,2.6) to[out=275, in=180] (3, 0.4) to[out=0, in=180] (4, 0.4);
\draw (1.4, 2.6) node[below right]{$\gamma$};
\draw[very thin] (1.72, 1.32) -- (0, 1.32);
\draw (-0.05, 1.32) node[left]{$\lambda^*(t)$};
\end{tikzpicture}
\caption{Given $t\in \R$, on the thick part of the line $\alpha-\beta=t$ we have existence and on the dashed part non-existence of a positive eigenfunction.}
\label{gamma}
\end{figure}

Partial results on $\sigma_{p, q}^+$ are scattered in the literature, since positive eigenfunctions turn out to be a useful tool for studying more general quasilinear problems. One always has
\[
]\lambda_1(p), +\infty[\  \times \ ]-\infty, \lambda_1(q)[ \ \cup\  ]-\infty,  \lambda_1(p)[ \ \times\  ]\lambda_1(q), +\infty[\ \subseteq \sigma_{p, q}^+,
\]
with uniqueness of positive eigenfunctions in $]-\infty, 0] \, \times\  ]\lambda_1(q), +\infty[$; cf. \cite[Proposition 2.2]{BT} as regards  existence and \cite[Lemma 2.2]{MaPa}, \cite[Theorem 1.1]{T} for uniqueness.\\
Moreover, a quite complete picture of $\sigma_{p, q}^+$ is available, according to whether the following conjecture holds true or not.
\begin{equation}\label{conj}
\begin{split}
&\text{Let $u_{1,p}$ and $u_{1,q}$ be first eigenfunctions for $A_p$ and $A_q$ respectively.}\\
&\text{If $p\neq q$ then $u_{1,p}$ and $u_{1,q}$ are linearly independent}.
\end{split}
\end{equation}
\begin{theorem}[\cite{BT}, Section 2]
Suppose $1<q<p<N$. Then
\[
\sigma_{p, q}\subseteq \{(\lambda_1(p), \lambda_1(q))\}\cup \{(\alpha,\beta):\,\alpha>\lambda_1(p)\}\cup\{(\alpha,\beta):\,\beta>\lambda_1(q)\}.
\]
If, moreover, \eqref{conj} is satisfied then there exists a continuous non-increasing function $\lambda^*:\R\to  [\lambda_1(q), +\infty[$
such that, letting $\gamma(t):=(\lambda^*(t)+t, \lambda^*(t))$,
\[
\sigma_{p, q}^+\setminus \gamma(\R)= 
\{(\alpha, \beta): \lambda_1(p)<\alpha, \beta<\lambda^*(\alpha-\beta)\} \cup\{(\alpha, \beta): \lambda_1(q)<\beta < \lambda^*(\alpha-\beta)\}
 \]
In addition, $\lambda^*(t)\equiv \lambda_1(q)$ for any $t$ big enough and $t\mapsto \lambda^*(t)+t$ is non-decreasing.
\end{theorem}
Therefore, the support of $\gamma$ plays the r\^{o}le of a threshold for the existence of positive eigenfunctions, as described in Figure \ref{gamma}. More precisely, any line
$$L_t:=\{(\alpha, \beta)\in \R^2: \alpha-\beta=t\}$$
intersects $\gamma$ at a unique point $P_t\equiv(\lambda^*(t)+t, \lambda^*(t))$ and $(\alpha, \beta)$ belongs to $L_t\cap (\sigma_{p,q}^+\setminus \gamma(\R))$ iff $(\alpha, \beta)$ lies in the part of $L_t$ {\em below}  $P_t$, under the constraint $\alpha>\lambda_1(p)$ or $\beta>\lambda_1(q)$.

Some natural properties of $\lambda^*$ are yet to be understood. As an example, is it true that $\lambda^*(t)+t\equiv \lambda_1(p)$ for $t$ sufficiently negative? Besides, while some partial results are known on the borderline case $\sigma_{p, q}^+\cap\gamma(\R)$, the picture looks not complete until today. Let us finally point out that, still in \cite{BT}, the case when Conjecture \eqref{conj} fails is discussed, providing a simpler description of $\sigma_{p, q}^+$. This is particularly meaningful for weighted eigenvalue problems, where it may occur that \eqref{conj} does not hold for some weights.
 
Regarding the set $\sigma_{p, q}\setminus \sigma_{p, q}^+$, or simply the existence of sign-changing solutions to \eqref{eig}, only partial results are available. 
\begin{theorem}\label{eigprobs}
\quad
\begin{enumerate}
\item Let $\alpha:=0$. Then \eqref{eig} admits a nodal eigenfunction iff $\beta>\lambda_2(q)$. See \cite[Theorem 3]{T0} and \cite[Theorem 1.2]{T}.
\item Assume that $\alpha<\lambda_1(p)$ and $\beta>\lambda_2(q)$. Then \eqref{eig} possesses a nodal eigenfunction. Cf. \cite[Theorem 4.3]{MaPa}.
\item Let $q:=2<p$. If $\alpha\notin \sigma_p$, $\beta\notin \sigma_2$, and 
\[
\alpha\notin [\lambda_{m(\beta)}(p), \lambda_{m(\beta)+1}(p)],
\]
where $m(\beta)$ is the number of $\lambda_k(2)\in \sigma_2$ fulfilling $\lambda_k(2)<\beta$, each counted with its geometric multiplicity, then $(\alpha, \beta)\in \sigma_{p, 2}$. See \cite[Theorem 4.2]{CDG}.
\end{enumerate}
\end{theorem}
Let us mention that the existence results of \cite{CDG,MaPa} actually deal with more general equations, as we will see later. Eigenvalue problems on the whole space are investigated in \cite{BZ0,BZ,HeLi}.
\section{Multiplicity results}
In this section we will discuss existence and multiplicity of solutions to the general quasilinear Dirichlet problem
\[
u\in W^{1,p}_0(\Omega),\quad A_{p, q}(u)=g(\cdot, u),
\]
where $g:\Omega\times\R\to\R$ satisfies Carath\'{e}odory's conditions. For the sake of simplicity, the case $g$ smooth and purely autonomous, i.e. $g(x, t)=\hat g(t)$ with $\hat g\in C^1(\R)$, is treated. Most of the results can be recasted as perturbed eigenvalue problems, which means
\begin{equation*}
\hat g(t):=\alpha |t|^{p-2}t+\beta |t|^{q-2}t +f(t), 
\end{equation*}
where the behavior of $f$ at $\pm \infty$ and/or at zero is typically negligible when compared with the leading term $\alpha |t|^{p-2}t+\beta |t|^{q-2}t$. To fix ideas, the sample equation is 
\begin{equation} \label{model}
A_{p, q}(u)=\alpha |u|^{p-2}u+\beta |u|^{q-2}u + \gamma |u|^{r-2}u
\end{equation}
for various $\alpha, \beta, \gamma \in \R$ and $r>1$, the last term representing the perturbation.

Let $1<q<p<N$ and let $f\in C^1(\R)$. We say that the problem 
\begin{equation}\label{pert}
u\in W^{1,p}_0(\Omega),\quad A_{p, q}(u)=\alpha |u|^{p-2}u+\beta |u|^{q-2}u + f(u)
\end{equation}
is {\em subcritical} if $f$ fulfills the growth rate
\[
|f(t)|\le C(1+|t|^{r-1})\quad\mbox{in}\quad\R,
\]
with appropriate $r\in]1,p^*[$ and {\em critical} when this inequality holds true for $r=p^*$ but not for smaller $r$'s. Moreover, $(p-1)$-sublinear signifies $r\leq p$, or, more generally,
\[
\limsup_{|t|\to+\infty}\frac{F(t)}{|t|^p}<+\infty,\quad\mbox{where}\quad F(t):=\int_0^t f(\tau)d\tau.
\]
Otherwise \eqref{pert} is called  $(p-1)$-superlinear. Finally, we say that Problem \eqref{pert} turns out to be {\em asymmetric} if 
$f$ (or $F$) exhibits a different asymptotic behavior at $\pm\infty$. 

Due to the expository nature of the present work, we won't examine here possible formulations of the Ambrosetti-Rabinowitz (briefly, (AR)) condition on the nonlinearities involved. This assumption, especially fruitful to achieve the Palais-Smale condition for $(p-1)$-superlinear problems, has been the object of many research papers in recent years. A deep discussion of various generalizations is made in \cite{LiYang} while \cite{CEM,MuPa,PaRa2} contain applications to $(p, q)$-Laplacian problems. We choose to somewhat oversimplify statements, substituting Condition (AR) with easier ones. So, the theorems given below are actually true under  less restrictive hypotheses than those taken on here.
\subsection{$(p-1)$-sublinear problems}
The present section collects some multiplicity results devoted to the sublinear case. Accordingly, the perturbation $f$ will fulfill
\[
-\infty<\lim_{|t|\to +\infty}\frac{f(t)}{|t|^{p-2}t}<+\infty.
\]
\subsubsection{Coercive setting}
This means that the energy functional stemming from \eqref{pert}, namely
\[
u\mapsto \frac{1}{p}\Vert u\Vert_p^p+\frac{1}{q}\Vert u\Vert_q^q-\frac{\alpha}{p}\vert u\vert_p^p-\frac{\beta}{q}\vert u\vert_q^q -\int_\Omega F(u(x))\, dx,\quad u\in W^{1,p}_0(\Omega),
\]
turns out to be coercive. As long as the behavior of $f(t)$ at infinity is negligible with respect to $|t|^{p-1}$, we are thus studying \eqref{pert} for $\alpha<\lambda_1(p)$.

Very recently, a multiplicity result has been proved in \cite{PRR} under coercivity conditions. The associated sample equation is \eqref{model} with $q=2$, $r\in \ ]2, p[$, $\gamma>0$.
\begin{theorem}[\cite{PRR}, Theorem 16]
Let $q:=2<p<N$ and let $f(0)=f'(0)=0$. If
\[
\lim_{|t|\to +\infty}\frac{f(t)}{|t|^{p-2}t}= 0, \quad \lim_{|t|\to +\infty} \frac{F(t)}{t^2}=+\infty
\]
then there exists $\eps>0$ such that, for every $(\alpha,\beta)\in\, ]\lambda_1(p)-\eps, \lambda_1(p)[\,\times(\, ]\lambda_2(2),+\infty[\, \setminus\sigma_2)$, Problem \eqref{pert} admits at least four nontrivial solutions. Moreover, one of them is positive and another negative.
\end{theorem}
The more general case $1<q<p$ is treated in \cite{MaPa}. It corresponds to $r\in [p,p^*[$ and $\gamma<0$  for the model equation \eqref{model}.
\begin{theorem}[\cite{MaPa}, Theorems 4.2 and 4.4]
Assume that $1<q<p<N$, $\alpha<\lambda_1(p)$, $\beta>\lambda_2(q)$, $f$ is subcritical, and
\[
\lim_{t\to 0}\frac{f(t)}{|t|^{q-2}t}=0,\quad\lim_{|t|\to +\infty}\frac{F(t)}{|t|^p}\le 0.
\]
Then \eqref{pert} possesses at least three solutions: $u_1\in{\rm int}(C_+)$, $u_2\in-{\rm int}(C_+)$, and $u_3\neq 0$. If, moreover, 
\[
\lim_{|t|\to +\infty}\frac{f(t)}{|t|^{p-2}t}>-\infty,
\]
then $u_3$ is nodal.
\end{theorem}
A concave perturbation is added in \cite{MePe}, obtaining the next multiplicity result.
\begin{theorem}[\cite{MePe}, Theorem 1.1]
Let $1<s<q<p<N$, let $\alpha<\lambda_1(p)$, $\beta>\lambda_1(q)$, and let $f$ be subcritical. If
\[
\lim_{t\to 0}\frac{f(t)}{|t|^{q-2}t}=0,\quad -\infty<\lim_{|t|\to +\infty} \frac{F(t)}{|t|^p}\le 0
\]
then there exists $\mu^*>0$ such that for every $\mu\in\,]0,\mu^*[$ the problem
\begin{equation*}
u\in W^{1,p}_0(\Omega),\quad A_{p, q}(u)=\alpha |u|^{p-2}u+\beta |u|^{q-2}u-\mu |u|^{s-2}u+f(u),
\end{equation*}
 has at least four nontrivial solutions, $u_1,u_2\in{\rm int}(C_+)$, $u_3,u_4\in-{\rm int}(C_+)$.
\end{theorem}
Under additional hypotheses on $F$, a fifth solution is found in the same paper.
\subsubsection{Resonant setting}
Roughly speaking, \eqref{pert} is called {\em resonant} provided $\alpha\in \sigma_p$ (or $\beta\in \sigma_q$) and $f$ is negligible at $\pm \infty$ with respect to $|t|^{p-1}$ (or $|t|^{q-1}$, respectively).

A problem resonant at higher eigenvalues has been addressed in \cite{PaRa2} for the $(p, 2)$-Laplacian, patterned after \eqref{model} with perturbation such that  $q<r<p$, $\gamma>0$.
\begin{theorem}[\cite{PaRa2}, Theorem 3.7]
Suppose $q:=2<r<p<N$, $\alpha:=\lambda_k(p)$ for some $k\ge 2$, and $\beta\in \ ]0,\lambda_1(2)[$. If $f$ satisfies 
\[
f(0)=f'(0)=0, \quad 0<\lim_{|t|\to +\infty}\frac{f(t)}{|t|^{r-2}t}<+\infty,\quad |f'(t)|\le C(1+|t|^{p-2})\;\;\forall\, t\in\R
\]
then \eqref{pert} admits at least three solutions: $u_1\in{\rm int}(C_+)$, $u_2\in-{\rm int}(C_+)$, and $u_3\neq 0$.
\end{theorem}
Concerning resonance at the first eigenvalue, we have the following
\begin{theorem}[\cite{MMP}, Theorems 3.9 and 4.5]
Let $1<s<q<r<p<N$, let $\alpha:=\lambda_1(p)$, $\beta:=0$, and let $f$ fulfill
\[
0<\lim_{t\to 0}\frac{f(t)}{|t|^{s-2}t}<+\infty.
\]
\begin{itemize}
\item {\em Resonance from the left.} Assuming
\[
 -\infty<\lim_{|t|\to +\infty}\frac{f(t)}{|t|^{r-2}t}<0
\]
yields three solutions: $u_1\in{\rm int}(C_+)$, $u_2\in-{\rm int}(C_+)$, and $u_3$ nodal.
\item {\em Resonance from the right.} On the other hand, the condition
\[
0<\lim_{|t|\to +\infty}\frac{f(t)}{|t|^{r-2}t}<+\infty,
\]
gives at least one nontrivial solution. 
\end{itemize}
\end{theorem}
Further papers on the same subject are \cite{BZ,MMP2,PaRa,PaWi2}; see also the references therein. In particular, \cite{BZ}  treats the case
$\Omega:=\R^N$, \cite{MMP2,PaRa} deal with asymmetric nonlinearities crossing an eigenvalue, while \cite{PaWi2} assumes $q:=2$.
\subsubsection{Problems neither coercive nor resonant}
Let us now come to the case where $f$ is $(p-1)$-sublinear but the energy functional associated with \eqref{pert} fails to be coercive and is indefinite. This occurs when, e.g.,
$$\lim_{|t|\to+\infty}\frac{F(t)}{|t|^p}>\lambda_1(p).$$
The following result has already been  mentioned in Section \ref{Sect3} for $f\equiv 0$.

\begin{theorem}[\cite{CDG}, Theorem 4.2]
Suppose $q:=2<p$, $\alpha\notin \sigma_p$, $\beta\notin \sigma_2$, $f(0)=f'(0)=0$, and
\[
\lim_{|t|\to +\infty}\frac{f'(t)}{|t|^{p-2}}=0.
\]
If, moreover,
\[
\alpha\notin [\lambda_{m(\beta)}(p), \lambda_{m(\beta)+1}(p)],
\]
where $m(\beta)$ is as in Theorem \ref{eigprobs}, then \eqref{pert} possesses  at least one nontrivial solution.
\end{theorem}
The sample equation for the next result is \eqref{model} with $q:=2$, $p<r$, $\gamma<0$. The setting is not purely resonant (meaning that $\alpha=\lambda_1(p)$ is not allowed in \eqref{pert}), but falls inside the so-called {\em near resonance} problems.
\begin{theorem}[\cite{PRR}, Theorems 22 and 28]
Let $q:=2<p<N$, let $f(0)=f'(0)=0$, and let $f$ be subcritical.
\begin{itemize}
\item If $(\alpha,\beta)\in \ ]\lambda_1(p), \lambda_2(p)[\ \times\ ]0, \lambda_1(2)[$ and
\[
0<\lim_{|t|\to +\infty}\frac{f(t)}{|t|^{r-2}t}<+\infty
\]
for appropriate $2<r<p$ then \eqref{pert} has a nontrivial solution.
\item Under the assumptions
\[
\lim_{|t|\to +\infty}\frac{F(t)}{|t|^{p}}< 0,\quad |f'(t)|\le C(1+|t|^{p-2})\;\;\forall\, t\in\R,
\]
there exists $\eps>0$ such that for every $(\alpha,\beta)\in\  ]\lambda_1(p), \lambda_1(p)+\eps[\ \times\ (\, ]\lambda_2(2),+\infty[ \,\setminus\sigma_2)$ Problem \eqref{pert} possesses at least four solutions: $u_1\in{\rm int}(C_+)$, $u_2\in-{\rm int}(C_+)$, and $u_3,u_4$ nodal.
\end{itemize}
\end{theorem}
Analogous results on the whole space, which causes further difficulties, are established in \cite{BZ,HeLi}. Finally, we refer to the survey paper \cite{MoWi} for non-variational problems involving the $(p, q)$-Laplace operator.
\subsection{$(p-1)$-superlinear problems}
This section contains some multiplicity results concerning the superlinear framework. So, the perturbation $f$ will fulfill
\[
\lim_{|t|\to +\infty}\frac{f(t)}{|t|^{p-2}t}=+\infty.
\]
\subsubsection{Subcritical setting}
Through the Nehari manifold approach, a {\em ground state} solution (namely a positive solution, which minimizes the associated energy functional) has been obtained  in \cite{FQ}.
\begin{theorem}[\cite{FQ}, Corollary 2.1]
Let $1<q<p<N$, let $\alpha:=\beta:=0$, and let $f$ be a subcritical function such that
\[
\lim_{t\to 0}\frac{f(t)}{|t|^{q-2}t}=0, \qquad t\mapsto \frac{f(t)}{t^{p-1}} \ \text{ is increasing on $\R^+$}.
\]
Then \eqref{pert} possesses a ground state solution.
\end{theorem}
When $q\le 2\le p$ we have the following
\begin{theorem}[\cite{MuPa}, Theorem 10]
Suppose $1<q\le 2\le p<N$, $\alpha<\lambda_1(p)$, $\beta:=0$. If $f$ satisfies
\[
-\infty<\lim_{t\to 0}\frac{f(t)}{|t|^{p-2}t}<0, \quad 0<\lim_{t\to +\infty}\frac{f(t)}{|t|^{r-2}t}<+\infty
\]
with appropriate $r\in \ ]p, p^*[$, then \eqref{pert} has three solutions: $u_1\in{\rm int}(C_+)$, $u_2\in-{\rm int}(C_+)$, and $u_3\neq 0$.
\end{theorem}
A more sophisticated result is proved in \cite{PaRa2}.
\begin{theorem}[\cite{PaRa2}, Theorem 4.12 and 4.11]
Let $q:=2<p<N$, let $f$ be a subcritical function, and let $(\alpha,\beta)\in\R^2$ satisfy 
\[
\alpha|c_+|^{p-2}c_++\beta c_++ f(c_+)\le 0\le \alpha |c_-|^{p-2}c_-+\beta c_-+ f(c_-)\quad \text{for some $c_-<0<c_+$}.
\]
Assume also that
\[
f(0)=f'(0)=0,\quad \lim_{|t|\to +\infty}\frac{f(t)}{|t|^{r-2}t}>0, \quad |f'(t)|\le C(1+|t|^{r-2})\;\;\forall\, t\in\R,
\]
where $p<r<p^*$. Then Problem \eqref{pert} admits: 
\begin{itemize}
\item Five solutions, $u_1,u_2\in{\rm int}(C_+)$, $u_3,u_4 \in-{\rm int}(C_+)$, $u_5\neq 0$, if $\beta\in\  ]\lambda_1(2), \lambda_2(2)[$.
\item Six solutions, two positive, two negative, and the other nodal, if $\beta\in\ ]\lambda_2(2),+\infty[\,\setminus\sigma_2$.
\end{itemize}
\end{theorem}
For additional results on the whole space, see \cite{BZ,CEM}.
\subsubsection{Critical setting}
In this framework, the most relevant term behaves as $|u|^{p^*-2}u$. Thus, we shall be concerned with equations of the type
\[
A_{p, q}(u)=|u|^{p^*-2}u+h(u), 
\]
where the perturbation $h$ is strictly subcritical. All the results we will present involve odd nonlinearities, so that any positive solution directly gives rise to a negative one. Under the hypothesis $1<r<q<p<N$, a first result on the sample problem
\begin{equation}\label{modelcrit}
u\in W^{1,p}_0(\Omega),\quad A_{p, q}(u)=|u|^{p^*-2}u+\gamma |u|^{r-2}u 
\end{equation}
has been obtained in \cite{LZ} and then generalized as follows.

\begin{theorem}[\cite{YY0}, Theorem 1.1]\label{YY0thm}
Let $1<r<q<p<N$ and let $f\in C^1(\R)$ be odd. If $f(t)t>0$ in $\R\setminus\{0\}$, there exists $s\in\ ]1,p[$ such that
\[
0<\lim_{|t|\to +\infty} \frac{f(t)}{|t|^{s-2}t}<+\infty,
\]
and $\gamma,\lambda>0$ are sufficiently small then the problem 
\begin{equation}
\label{yy}
u\in W^{1,p}_0(\Omega),\quad A_{p, q}(u)=|u|^{p^*-2}u+\gamma|u|^{r-2}u+\lambda f(u)
\end{equation}
possesses infinitely many solutions.
\end{theorem}
Superlinear perturbations of the purely critical equation are treated in \cite{YY}.
\begin{theorem}[\cite{YY}, Theorems 1.1 and 1.2]
Suppose $1<q<p<r<p^*$ and $p<N$. Then:
\begin{itemize}
\item \eqref{modelcrit} admits a nontrivial solution provided $\gamma>0$ is large enough.
\item Let  ${\rm cat}(\Omega)$ denote the Ljusternik--Schnirelmann category of $\Omega$ in itself. If
\[
N>p^2, \quad p^*-r< \frac{q}{p-1}<1^*,
\]
then \eqref{modelcrit} has at least ${\rm cat}(\Omega)$ positive solutions for every sufficiently small $\gamma>0$. 
\end{itemize}
\end{theorem}
Let us note that Theorem \ref{YY0thm} actually extends to more general equations of the form \eqref{yy}, still for $r\in \ ]p, p^*[$, cf. \cite[Theorem 4.3--4.4]{YY0}.

Finally, the borderline case of an eigenvalue problem with critical nonlinearity is investigated by following result.
\begin{theorem}[\cite{CMP}, Theorem 1.3]
Let $1<q<p<\min\{N, q^*\}$ and let $\beta\in \R$. Then, for every $\alpha>0$ large enough, the problem
\[
u\in W^{1,p}_0(\Omega),\quad A_{p, q}(u)=\alpha |u|^{p-2}u+\beta |u|^{q-2}u+|u|^{p^*-2}u
\]
possesses a nontrivial solution, which is strictly positive provided $\beta<\lambda_1(q)$.
\end{theorem}

Then paper \cite{CMP} contains further existence results concerning the situation
$$1<q<p<q^*<N.$$
It should be noted that, even for the sample problem \eqref{modelcrit}, the intermediate case $q<r<p$ is, as far as we know, still open. Finally, the existence of {\em ground states} for critical equations on the whole space is studied in \cite{F}.
\subsection{Asymmetric nonlinearities}
We now discuss quasilinear Dirichlet problems with reactions having different asymptotic behaviors at  $-\infty$ and $+\infty$. The sample equation stems from the Fu\v cik spectrum theory, and is of the type
\begin{equation}\label{fucik}
A_{p, q}(u)=\alpha_+ u_+^{p-1}-\alpha_- u_-^{p-1} +\beta |u|^{q-2}u+f(u), \qquad\alpha_\pm, \beta\in \R,
\end{equation}
where $t_\pm:=\max\{0, \pm t\}$ while $f$ exhibits suitable (possibly asymmetric as well) growth rates at zero and $\pm\infty$. For the sake of simplicity, we singled out the Fu\v cik structure only on the higher order term of the reaction with respect to which resonance can occur, but, as before, most of theorems can be recasted in a more general setting. 

A first result in this framework is obtained provided $1<q<p:=2$, both for subcritical and critical nonlinearities.
\begin{theorem}[\cite{YangBai}, Theorems 1.1 and 1.2]
Suppose $1<q<s<p:=2<r\le 2^*$, $\alpha\in\ ]\lambda_1(2),+\infty[\,\setminus\sigma_2$, and $\lambda>0$. Then the problem
\[
u\in W^{1,2}_0(\Omega),\quad A_{2, q}(u)=\alpha u-\lambda |u|^{s-2}u+\mu u_+^{r-1}
\]
admits at least three nontrivial solutions if $\mu>0$ is sufficiently small and either $r<2^*$, $N\ge 3$ or $r=2^*$, $N\ge 4$.
\end{theorem}
The next result allows a full resonance (from the left) at $+\infty$ with respect to $\lambda_1(p)$.
\begin{theorem}[\cite{PaRa3}, Theorems 3.4 and  4.3]
Let $1<q:=2<r<p<N$, let $\alpha_-\leq \lambda_1(p)<\alpha_+$, and let $\beta\in\ ]\lambda_2(2),+\infty[\,\setminus\sigma_2$. If $f(0)=f'(0)=0$,
\[
-\infty<\lim_{t\to -\infty}\frac{f(t)}{|t|^{r-2}t}<0<\lim_{t\to +\infty}\frac{f(t)}{t^{r-1}}<+\infty,\quad\mbox{and}\quad
|f'(t)|\le C(1+|t|^{p-2})\;\;\forall\, t\in\R
\]
then \eqref{fucik} possesses at least two nontrivial solutions, one of which is negative. A third solution exists once $\alpha_-\neq \lambda_1(p)$.
\end{theorem}
When the perturbation in \eqref{fucik} contains a parametric concave term we have the following
\begin{theorem}[\cite{MMP2}, Theorem 4.2]
Suppose $1<s<q<p<N$ and $\alpha_-\le \lambda_1(p)\le \alpha_+$. If, moreover,
\[
0<\lim_{t\to 0}\frac{f(t)}{|t|^{q-2}t}<\lambda_1(q),\quad -\infty<\lim_{t\to -\infty}\frac{f(t)}{|t|^{p-2}t}<0<\lim_{t\to +\infty}\frac{f(t)}{t^{p-1}}<+\infty,
\]
then there exists $\lambda^*>0$ such that, for every $\lambda\in \ ]0, \lambda^*[$, the  problem
\begin{equation*}\label{par}
u\in W^{1,p}_0(\Omega),\quad A_{p, q}(u)=\alpha_+ u_+^{p-1}-\alpha_- u_-^{p-1} +\lambda |u|^{s-2}u+f(u),
\end{equation*}
admit at least four solutions, $u_1,u_2\in{\rm int}(C_+)$, $u_3\in-{\rm int}(C_+)$, and $u_4$ nodal.
\end{theorem}
Let us finally point out that infinitely many solutions are also obtained in \cite{MMP2} under a symmetry condition near zero. For further results concerning asymmetric nonlinearities, see \cite{PaRa4,PeiZhang}.
\section*{Acknowledgement}
Work performed under the auspices of GNAMPA of INDAM.
\end{document}